# Optimal Thrust Level for Orbit Insertion

Max CERF [*]

Airbus Defence and Space, France

## Abstract

The minimum-fuel orbital transfer is analyzed in the case of a launcher upper stage using a constantly thrusting engine. The thrust level is assumed to be constant and its value is optimized together with the thrust direction. A closed-loop solution for the thrust direction is derived from the extremal analysis for a planar orbital transfer. The optimal control problem reduces to two unknowns, namely the thrust level and the final time. Guessing and propagating the costates is no longer necessary and the optimal trajectory is easily found from a rough initialization. On the other hand the initial costates are assessed analytically from the initial conditions and they can be used as initial guess for transfers at different thrust levels. The method is exemplified on a launcher upper stage targeting a geostationary transfer orbit.

## 1. Introduction

Mass minimization is a major concern for the design of launch vehicles. The fuel required to reach the targeted orbit depends on both the thrust level and the thrust orientation along the trajectory. Finding the minimum-fuel trajectory is an optimal control problem that can be addressed by the Pontryaguin Maximum Principle (PMP). Due to its utmost importance, this problem has received a considerable attention from academics and industrials since the beginning of the space age.

In a few cases, when the dynamical model is sufficiently simple, an analytical solution may be derived from the PMP necessary conditions. Among such well-known solutions for minimum-fuel trajectory problems, we can mention the following ones.

- The Goddard's problem [1] for a vertical launched rocket.
- The flat Earth model with constant gravity and constant acceleration [2,3]
- The Edelbaum's model for a low thrust transfer between circular orbits [4]

[*] AIRBUS Defence & Space, 78130 Les Mureaux, France
max.cerf@astrium.eads.net

The flat Earth model and the constant acceleration model do not have the sufficient representativeness to correctly assess the launcher optimal performance [5]. The problem must be formulated in a central gravity field and considering the actual engine thrust level. Depending on the launcher and the mission specifications, the thrust level may be either prescribed (in that case the problem is equivalent to a minimum-time problem) or freely optimized between some bounds. With these dynamical assumptions, there is no analytical solution and numerical methods must be used either direct or indirect [6,7].

Applying the PMP to the minimum-fuel problem yields in the regular case the optimal thrust direction aligned with the velocity costate, whereas the thrust level (if optimized) is driven by the sign of a switching function. Such switches induce additional issues : their number is a priori unknown and they must be detected accurately within the simulation process in order to keep a good differentiability of the problem [8]. Furthermore there may exist singular arcs along which the PMP first order conditions no longer define the optimal control. Such singular solutions require further theoretical analysis and specific solution methods [9].
Even when only the thrust direction is optimized (regular case, prescribed thrust level), finding the unknown initial costate proves numerically challenging, and it may discourage from using an indirect method.

This paper addresses the minimum-fuel orbital transfer in the particular case of an engine constantly thrusting at the same thrust level. Although this assumption may seem restrictive, this case is of great practical importance since rocket engines are generally designed for a reference thrust level. The targeted application is the flight of a launcher upper stage. The initial conditions are prescribed resulting from the previous stage flight. The final conditions are defined in terms of orbital parameters.
When the engine constant thrust level is considered as a free optimization parameter (as is the case in preliminary design studies) an additional optimality condition has to be written with the PMP equations. This condition can be exploited in the planar case to derive a closed-loop solution for the thrust optimal direction. Guessing and propagating the costates is no longer necessary and the minimum-fuel trajectory problem is reduced to a nonlinear system of 2 equations (targeted apogee and perigee) with 2 unknowns (thrust level and final time). This problem is easily solved from a rough initial solution.
On the other hand the initial costates corresponding to the optimal trajectory are derived analytically from the initial conditions. These costates can be used as initial guess for instances of the same minimum-fuel problem when the thrust level is no longer a free parameter.

The text is organized as follows. In section §2 the optimal control problem is formulated and the extremal conditions are analyzed. A closed-loop control law is derived in the planar case and the solution method is presented. In section §3 the method is applied to a representative example of a launcher targeting a geostationary transfer orbit. A sensitivity analysis on the thrust level illustrates how the analytical costates can be used as starting point to solve the minimum-fuel problem at non-optimal thrust level. The extension to low-thrust transfers is also discussed.

## 2. Problem Formulation and Analysis

This section formulates the Optimal Control Problem (OCP) under consideration. The problem is analyzed by applying the Pontryaguin Maximum Principle (PMP) and a closed-loop control is derived from the first order necessary conditions in the planar case.

### 2.1 Optimal control problem

The problem consists in finding the minimal-fuel trajectory to go from given injection conditions to a targeted orbit with a constantly thrusting engine. The Earth is modeled as a sphere, with its center at the origin of an inertial frame. The vehicle is considered as a material weighting point with position $\vec{r}(t)$, velocity $\vec{v}(t)$, mass $m(t)$ submitted to the Earth acceleration gravity denoted $\vec{g}(\vec{r})$ and to the engine thrust. The thrust level T is constant with a burned propellant exhaust velocity equal to $v_e$. The engine is ignited at the initial date $t_0$ and it cannot be turned off before the orbit insertion at $t_f$. The thrust direction can be chosen freely and it is orientated along the unit vector $\vec{u}(t)$.

Applying the fundamental dynamics principle in the Earth-centered inertial frame yields the motion equations.

$$\begin{cases} \dot{\vec{r}} = \vec{v} \\ \dot{\vec{v}} = \vec{g} + \dfrac{T}{m}\vec{u} \\ \dot{m} = -\dfrac{T}{v_e} \end{cases} \tag{1}$$

The dependencies on time (for $\vec{r}$, $\vec{v}$, $m$, $\vec{u}$) and on position (for $\vec{g}$) have been omitted for conciseness.

In order to formulate an optimal control problem, we consider as state variables $\vec{r}(t)$, $\vec{v}(t)$ and $m(t)$. The control variables are the thrust direction $\vec{u}(t)$, the thrust level T and the final time $t_f$. The thrust level T can take any positive value (no upper bound).

The initial state at the engine ignition results from the ascent trajectory flown by the launcher lower stages. This initial state is completely prescribed. The final state is constrained by the targeted orbit defined by the apogee and perigee altitudes denoted respectively $h_A$ and $h_P$. The apogee and perigee altitudes actually achieved at the final date are denoted respectively $\psi_A$ and $\psi_P$ and they depend on the final position $\vec{r}(t_f)$ and velocity $\vec{v}(t_f)$

The optimal control problem is formulated under the Mayer form with a final cost.

$$\min_{\vec{u}(t), T, t_f} J = -m(t_f) \quad \text{s.t.} \quad \begin{cases} \dot{\vec{r}} = \vec{v} \\ \dot{\vec{v}} = \vec{g} + \dfrac{T}{m}\vec{u} \\ \dot{m} = -\dfrac{T}{v_e} \end{cases} \text{with} \begin{cases} \vec{r}(t_0) = \vec{r}_0 \\ \vec{v}(t_0) = \vec{v}_0 \\ m(t_0) = m_0 \end{cases} \text{fixed initial state} \\ \begin{cases} \psi_A[\vec{r}(t_f), \vec{v}(t_f)] = h_A \\ \psi_P[\vec{r}(t_f), \vec{v}(t_f)] = h_P \end{cases} \text{constrained final state} \tag{2}$$

## 2.2 Extremal Analysis

The optimal trajectory is sought by applying the Pontryaguin Maximum Principle (PMP) [10,11]. For that purpose we introduce the costate vectors $\vec{p}_r(t)$, $\vec{p}_v(t)$, $p_m(t)$ respectively associated to the position, the velocity, the mass. These costate vectors do not vanish identically on any interval of $[t_0, t_f]$ and they are defined up to a non-positive scalar multiplier $p_0$ which can be chosen freely for a normalization purpose.

With these notations, the Hamiltonian for the OCP Eq. (2) is

$$H = \vec{p}_r^T \dot{\vec{r}} + \vec{p}_v^T \dot{\vec{v}} + p_m \dot{m} = \vec{p}_r^T \vec{v} + \vec{p}_v^T \left( \vec{g} + \frac{T}{m} \vec{u} \right) + p_m \left( -\frac{T}{v_e} \right)$$
$$= \vec{p}_r^T \vec{v} + \vec{p}_v^T \vec{g} + T \left( \frac{\vec{p}_v^T \vec{u}}{m} - \frac{p_m}{v_e} \right) \qquad (3)$$

The PMP provides the following first order necessary conditions on $\vec{u}$, T and $t_f$ to be an optimal control.

- The Hamiltonian maximization condition with respect to the control $\vec{u}(t)$.

$$\max_{\vec{u}} H \qquad (4)$$

  This condition leads to a thrust direction aligned with the velocity costate.

$$\vec{u} = \frac{\vec{p}_v}{p_v} \quad \text{with} \quad p_v = \|\vec{p}_v\| \qquad (5)$$

- The costate differential equations.

$$\begin{cases} \dot{\vec{p}}_r = -\frac{\partial H}{\partial \vec{r}} = -\frac{\partial \vec{g}}{\partial \vec{r}} \vec{p}_v \\ \dot{\vec{p}}_v = -\frac{\partial H}{\partial \vec{v}} = -\vec{p}_r \\ \dot{p}_m = -\frac{\partial H}{\partial m} = \frac{T}{m^2} \vec{p}_v^T \vec{u} = \frac{T}{m^2} p_v \quad \text{using} \quad \vec{u} = \frac{\vec{p}_v}{p_v} \end{cases} \qquad (6)$$

- The transversality conditions on the final costate, derived from the final constraints $\psi_A$, $\psi_P$ with the multipliers $\nu_A$, $\nu_P$ and from the final cost $-m(t_f)$ with the non-positive multiplier $p_0$.

$$\begin{cases} \vec{p}_r(t_f) = \nu_A \frac{d\psi_A}{d\vec{r}}(t_f) + \nu_P \frac{d\psi_P}{d\vec{r}}(t_f) \\ \vec{p}_v(t_f) = \nu_A \frac{d\psi_A}{d\vec{v}}(t_f) + \nu_P \frac{d\psi_P}{d\vec{v}}(t_f) \\ p_m(t_f) = -p_0 \end{cases} \qquad (7)$$

- The transversality condition on the free final date $t_f$.

$$H(t_f) = 0 \qquad (8)$$

  The problem is furthermore autonomous, so that the Hamiltonian is constant along any extremal.

$$H(t) = 0 \quad, \quad \forall t \in [t_0, t_f] \tag{9}$$

- The first order condition on the scalar parameter T. The proof for this condition can be found in [12,13].

$$\int_{t_0}^{t_f} \frac{\partial H}{\partial T} dt = 0 \quad \Rightarrow \quad \int_{t_0}^{t_f} \left( \frac{p_v}{m} - \frac{p_m}{v_e} \right) dt = 0 \tag{10}$$

We now search for a closed-loop control law $\vec{u}(t)$ starting from Eq. (10). The dependencies on the time are omitted for conciseness, but they are explicitly written whenever necessary for a better comprehension.

By Bellman's principle of optimality [11,13] and since the problem is autonomous, the condition Eq. (10) must in fact hold from any date t along an optimal trajectory, and not only from the initial date $t_0$.

$$\int_{t}^{t_f} \left( \frac{p_v}{m} - \frac{p_m}{v_e} \right) dt = 0 \quad, \quad \forall t \in [t_0, t_f] \tag{11}$$

This implies that the integrand is constantly null along an optimal trajectory.

$$\frac{p_v}{m} - \frac{p_m}{v_e} = 0 \quad, \quad \forall t \in [t_0, t_f] \tag{12}$$

Using the third equation in Eq. (6) to replace $p_v$, we can integrate and find an expression for $p_m$.

$$\frac{m \dot{p}_m}{T} - \frac{p_m}{v_e} = 0 \quad \Rightarrow \quad m \dot{p}_m + \dot{m} p_m = 0 \qquad \text{using } \dot{m} = -\frac{T}{v_e} \text{ from Eq.(2)}$$
$$\Rightarrow \quad \frac{d}{dt}(p_m m) = 0 \tag{13}$$
$$\Rightarrow \quad p_m(t) m(t) = C^{te} = p_m(t_f) m(t_f) \quad, \quad \forall t \in [t_0, t_f]$$

The value of $p_m(t_f)$ is given by the transversality condition Eq. (7).

$$p_m(t) m(t) = -p_0 m(t_f) \quad, \quad \forall t \in [t_0, t_f] \tag{14}$$

Replacing $p_m$ in Eq. (12) we find that the modulus of the velocity costate is constant.

$$p_v(t) = -\frac{p_0 m(t_f)}{v_e} \quad, \quad \forall t \in [t_0, t_f] \tag{15}$$

The constant $p_v$ can not be null else from Eq. (6) $\vec{p}_r$, $\vec{p}_v$, $p_m$ would vanish simultaneously, which would be in contradiction with the PMP. $p_0$ is therefore strictly negative and there exist no abnormal extremal. We can choose the usual normalization factor $p_0 = -1$.

This constant modulus property can be used to find the optimal thrust direction $\vec{u}(t)$ in the case of a planar trajectory. For that purpose, we come back to the costate equations Eq. (6) which yield a second order differential equation for $\vec{p}_v$.

$$\begin{cases} \dot{\vec{p}}_r = -\frac{\partial \vec{g}}{\partial \vec{r}} \vec{p}_v \\ \dot{\vec{p}}_v = -\vec{p}_r \end{cases} \quad \Rightarrow \quad \ddot{\vec{p}}_v = \frac{\partial \vec{g}}{\partial \vec{r}} \vec{p}_v \tag{16}$$

The vector $\vec{p}_v(t)$ has a constant modulus $p_v$ and it is rotating with some unknown angular rate $\omega(t)$ in the Galilean frame. Its derivative $\dot{\vec{p}}_v(t)$ is a vector normal to $\vec{p}_v(t)$ with modulus $p_v\omega(t)$. Denoting $\vec{n}$ the vector normal to the thrust direction $\vec{u}$ in the trajectory plane with the positive orientation, we have

$$\vec{p}_v = p_v\vec{u} \Rightarrow \dot{\vec{p}}_v = p_v\omega\vec{n} \Rightarrow \ddot{\vec{p}}_v = p_v(\dot{\omega}\vec{n} - \omega^2\vec{u}) \quad \text{using} \quad \begin{cases} \dot{\vec{u}} = \omega\vec{n} \\ \dot{\vec{n}} = -\omega\vec{u} \end{cases} \quad (17)$$

The Figure 1 depicts the vectors and angles used for the subsequent calculations. $\varphi$ is the longitude measured from the x axis in the Galilean frame (O,x,y). $\theta$ and $\gamma$ are respectively the local pitch angle and the flight path angle measured positively from the local horizontal. $\omega$ is the angular rate of $\vec{u}$ in the Galilean frame (O,x,y).

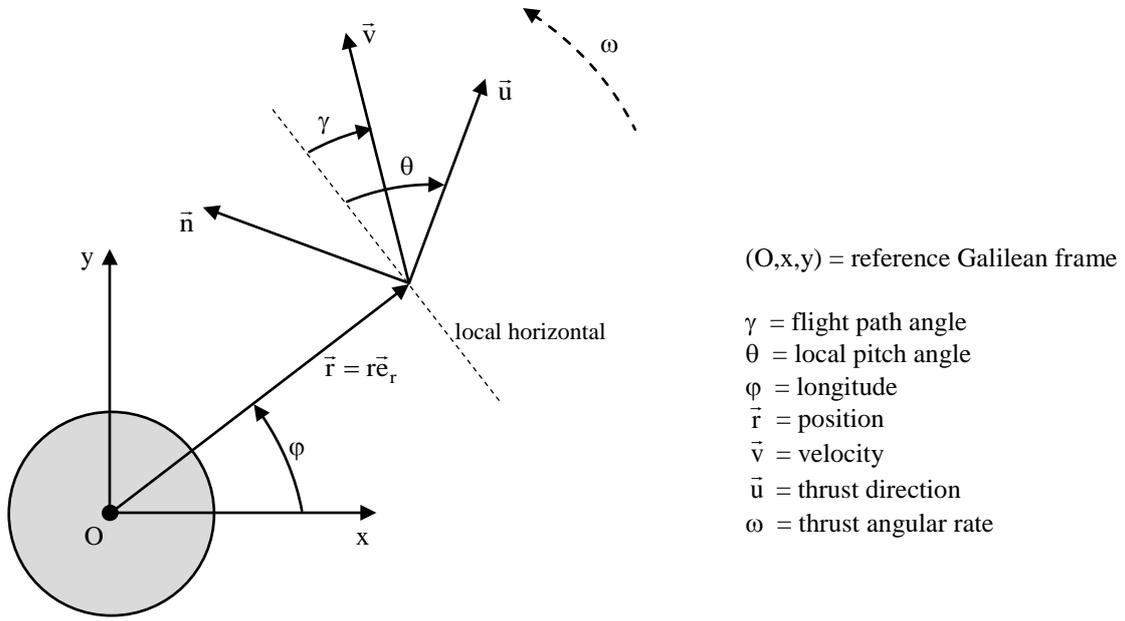

(O,x,y) = reference Galilean frame

$\gamma$ = flight path angle
$\theta$ = local pitch angle
$\varphi$ = longitude
$\vec{r}$ = position
$\vec{v}$ = velocity
$\vec{u}$ = thrust direction
$\omega$ = thrust angular rate

Figure 1 : Vectors and angles

Assuming a central gravity field, the gradient gravity matrix denoted G has an explicit expression depending on the position $\vec{r}$ [14,15] (this formula is obtained by deriving the gravity field in cartesian coordinates).

$$\vec{g}(\vec{r}) = -\frac{\mu}{r^3}\vec{r} \Rightarrow G = \frac{\partial\vec{g}}{\partial\vec{r}} = \frac{\mu}{r^5}(3\vec{r}\vec{r}^T - r^2 I) = \frac{\mu}{r^3}(3\vec{e}_r\vec{e}_r^T - I) \quad (18)$$

where I is the identity matrix and $\vec{e}_r$ is the unit radial vector : $\vec{r} = r\vec{e}_r$.

Replacing in Eq. (16) $\ddot{\vec{p}}_v$ given by Eq. (17) and G given by Eq. (18), and simplifying by $p_v$ which is not null

$$\dot{\omega}\vec{n} - \omega^2\vec{u} = \frac{\mu}{r^3}(3\vec{e}_r\vec{e}_r^T - I)\vec{u} = \frac{3\mu}{r^3}(\vec{e}_r^T\vec{u})\vec{e}_r - \frac{\mu}{r^3}\vec{u} \quad (19)$$

The projection of $\vec{e}_r$ on $(\vec{u},\vec{n})$ is (see Figure 1) : $\vec{e}_r = \vec{u}\sin\theta - \vec{n}\cos\theta$, so that the equation becomes

$$\vec{u}\left[\omega^2 + \frac{\mu}{r^3}(3\sin^2\theta - 1)\right] - \vec{n}\left[\dot{\omega} + \frac{3\mu}{r^3}\sin\theta\cos\theta\right] = \vec{0} \quad (20)$$

We obtain two relationships linking the local pitch angle θ and the thrust angular rate ω in the Galilean frame.

$$\begin{cases} \omega^2 = \dfrac{\mu}{r^3}(1-3\sin^2\theta) \\ \dot{\omega} = -\dfrac{3\mu}{r^3}\sin\theta\cos\theta \end{cases} \tag{21}$$

The first equation bounds the pitch angle since (1 − 3sin²θ) must be positive.

$$\sin^2\theta \leq \frac{1}{3} \quad \Rightarrow \quad \begin{cases} -35.26\deg \leq \theta \leq 35.26\deg \\ \text{or} \\ 144.73\deg \leq \theta \leq 215.26\deg \end{cases} \tag{22}$$

Another expression for the angular rate ω can be derived from the Hamiltonian nullity Eq. (9).

$$H = H_0 + TH_T = 0 \quad \text{with} \quad \begin{cases} H_0 = \vec{p}_r^T\vec{v} + \vec{p}_v^T\vec{g} \\ H_T = \dfrac{\vec{p}_v^T\vec{u}}{m} - \dfrac{p_m}{v_e} = \dfrac{p_v}{m} - \dfrac{p_m}{v_e} \end{cases} \tag{23}$$

From Eq. (12), we have : $H_T = 0$, so that $H_0$ must also vanish.

$$\vec{p}_r^T\vec{v} + \vec{p}_v^T\vec{g} = 0 \quad \text{with} \quad \begin{cases} \vec{p}_v = p_v\vec{u} \\ \vec{p}_r = -\dot{\vec{p}}_v = -\omega p_v\vec{n} \\ \vec{v} = v(\vec{u}\cos(\theta-\gamma)+\vec{n}\sin(\theta-\gamma)) \\ \vec{g} = -g\vec{e}_r = -g(\vec{u}\sin\theta - \vec{n}\cos\theta) \end{cases} \text{and} \quad g = \dfrac{\mu}{r^2} \tag{24}$$

This yields after simplification by $p_v$.

$$\omega v \sin(\theta-\gamma) - \frac{\mu}{r^2}\sin\theta = 0 \tag{25}$$

Replacing ω given by Eq. (21), and assuming that θ is within the bounds Eq. (22) we obtain a relationship between the thrust direction θ and the kinematic conditions (r, v, γ).

$$\begin{cases} \omega = \pm\sqrt{\dfrac{\mu}{r^3}(1-3\sin^2\theta)} \\ \omega v \sin(\theta-\gamma) = \dfrac{\mu}{r^2}\sin\theta \end{cases} \Rightarrow \sqrt{1-3\sin^2\theta}\, v\sin(\theta-\gamma) = \pm\sqrt{\dfrac{\mu}{r}}\sin\theta \tag{26}$$

Denoting $v_c$ the circular velocity at the radius r this equation can be written

$$\sin(\theta-\gamma) = \pm\frac{v_c}{v}\frac{\sin\theta}{\sqrt{1-3\sin^2\theta}} \quad \text{with} \quad v_c = \sqrt{\frac{\mu}{r}} \tag{27}$$

This equation is a necessary condition that must be fulfilled at any date by the optimal thrust direction θ.
It defines a closed-loop control since it depends only on the current kinematic conditions (r, v, γ).

The solution θ can exist only if the second member is comprised between −1 and +1.

$$\frac{v_c}{v}\frac{|\sin\theta|}{\sqrt{1-3\sin^2\theta}} \leq 1 \quad \Rightarrow \quad \sin^2\theta \leq \frac{1}{3+\left(\dfrac{v_c}{v}\right)^2} \tag{28}$$

This defines tighter bounds than Eq. (22) on the pitch angle θ.

$$\begin{cases} -\theta_m \leq \theta \leq \theta_m \\ \text{or} \\ \pi - \theta_m \leq \theta \leq \pi + \theta_m \end{cases} \quad \text{with} \quad \theta_m = \text{Arcsin} \frac{1}{\sqrt{3 + \left(\frac{v_c}{v}\right)^2}} \quad , \quad 0 < \theta_m < \frac{\pi}{2} \tag{29}$$

We observe that if θ is a solution of Eq. (27), then θ+π is also a solution. We can restrict ourselves to searching the solution in the interval $[-\theta_m, +\theta_m]$ and deduces afterwards all the possible solutions of Eq. (27).

Eq. (27) has no analytical solution and a numerical procedure is needed. For that purpose, γ is expressed as a function of θ by taking the Arc sin of both members.

$$\gamma = \theta \pm \text{Arcsin}\left(\frac{v_c}{v} \frac{\sin\theta}{\sqrt{1 - 3\sin^2\theta}}\right) + k\pi \quad , \quad k \text{ integer} \tag{30}$$

From the previous remark the solution θ is defined modulo π. It is thus sufficient to solve Eq. (30) for k=0.
We have finally two equations to solve in the interval $[-\theta_m, +\theta_m]$.

$$\gamma = F^+(\theta) \quad \text{or} \quad \gamma = F^-(\theta) \tag{31}$$

where $F^+(\theta)$ and $F^-(\theta)$ design the second member respectively with a plus or minus sign.

By inspection, we observe that these fonctions $F^+(\theta)$ and $F^-(\theta)$ are monotonous in the interval $[-\theta_m, +\theta_m]$, the ratio $v_c/v$ being in practice comprised between 0.7 (escape velocity) and 1.5 (low initial velocity). The Figure 2 present plots of the functions $F^+$ (left) and $F^-$ (right) for three values of this velocity ratio.

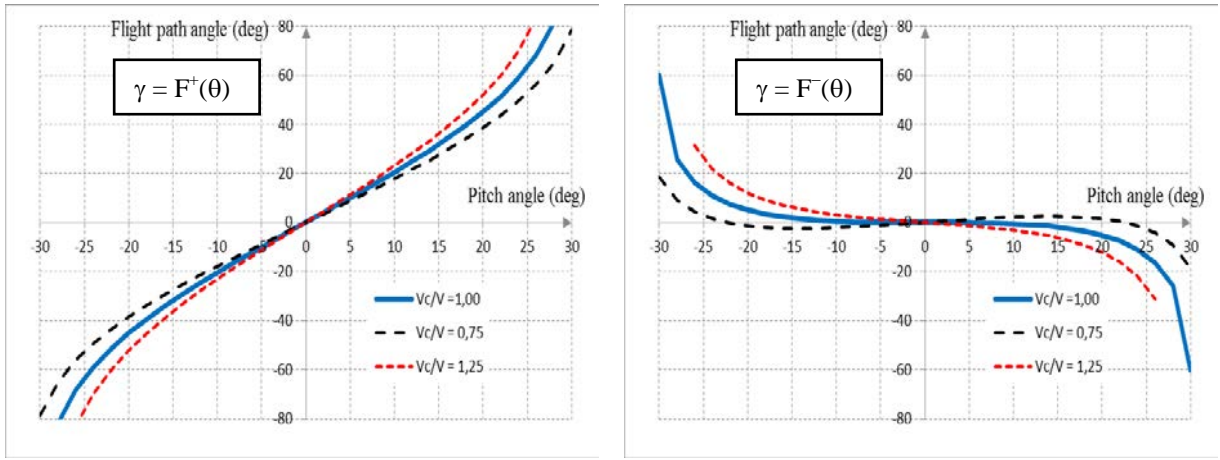

Figure 2 : Functions $F^+$ (left) and $F^-$ (right)

For a given value γ of the flight path angle, it appears that the first equation $F^+(\theta) = \gamma$ has always a unique solution denoted $\theta^+$. This is not the case for the second equation $F^-(\theta) = \gamma$. Furthermore this first solution $\theta^+$ is more physical since it follows the flight path angle evolution. We will therefore consider only this solution for the practical application.

The conclusion of this analysis of the minimum-fuel transfer in the planar case can now be drawn.

At any date t with given kinematic conditions (r, v, γ), there exist either 2 or 4 candidate thrust directions satisfying the PMP necessary conditions.

- A first direction $\theta^+$ is obtained by solving the Eq. (32). This solution always exists and it corresponds to a thrust direction in the vicinity of the velocity direction. The opposite direction $\theta^+ + \pi$ is also a solution.

$$\gamma = \theta + \text{Arcsin}\left(\frac{v_c}{v} \frac{\sin\theta}{\sqrt{1-3\sin^2\theta}}\right) \quad \text{with} \quad -\theta_m \leq \theta \leq \theta_m \quad , \quad \theta_m = \text{Arcsin}\frac{1}{\sqrt{3+\left(\frac{v_c}{v}\right)^2}} \quad , \quad 0 < \theta_m < \frac{\pi}{2} \quad (32)$$

- A second direction $\theta^-$ may exist considering the equation with a minus sign. This solution although satisfying the PMP is not physically meaningful.

For a launcher trajectory, the only meaningful solution is thus the first one $\theta^+$ which contributes to a velocity increase. The other candidate solutions can be discarded. This is not the case for low thrust transfers which may comprise arcs with thrusting opposite to the velocity. The correct solution should then be selected from other considerations regarding the whole trajectory and the endpoint conditions.

### 2.3 Solution Method

The unknowns of OCP Eq. (2) are the thrust direction $\vec{u}(t)$, the thrust level T and the final time $t_f$. For a planar trajectory, Eq. (32) allows finding the thrust pitch angle θ directly at any date from the current kinematic conditions (r, v, γ). This equation has a unique solution in the bounds considered. It is solved numerically "on-line" by a dichotomy method while propagating the motion differential equations Eq. (1).

The thrust direction $\vec{u}(t)$ is therefore no longer part of the problem unknowns. The OCP Eq. (2) reduces to a nonlinear system of two equations (final apogee and perigee) whose unknowns are the thrust level T and the final time $t_f$. This system is easily solved using a Newton-like method, without any particular care regarding the initial guess.

The initial guess may for example be based on the total velocity impulse required by the transfer. The required velocity impulse is the sum of the actual velocity increment and the velocity losses incurred along the trajectory.

$$\Delta V = v_e \ln\frac{m_0}{m_f} = v_f - v_0 + \Delta V_{loss} \quad (33)$$

The actual velocity increment $v_f - v_0$ can be estimated by assuming that the injection will occur at the perigee of the targeted orbit. Ignoring the velocity losses $\Delta V_{loss}$ (unknown a priori) yields an initial guess for $m_f$.

$$m_f \approx m_0 e^{-\frac{v_f - v_0}{v_e}} \quad (34)$$

This value of $m_f$ can in turn be used to build an initial guess for T based on a typical value $a_f$ of the final acceleration, between 10 m/s² and 20 m/s².

$$T \approx m_f a_f \tag{35}$$

The final time $t_f$ is then deduced from the mass consumption.

$$t_f = \frac{(m_0 - m_f)v_e}{T} \approx \frac{v_e}{a_f}\left(e^{\frac{v_f - v_0}{v_e}} - 1\right) \tag{36}$$

The convergence on the optimal solution is readily achieved from these rough estimates for T and $t_f$.

Once the nonlinear system is solved, we obtain the maximum final mass $m_f^*$, the optimal thrust level $T^*$, and the optimal command law $\theta^*(t)$ along the trajectory. The initial costates $\vec{p}_r(t_0)$, $\vec{p}_v(t_0)$, $p_m(t_0)$ of the optimal control problem Eq. (2) can then be retrieved as follows.

The normalization factor is set to the usual value $p_0 = -1$ (the costate vector is defined up to this factor).
The initial mass costate is assessed from Eq. (14).

$$p_m(t_0) = \frac{m_f^*}{m_0} \tag{37}$$

The initial position and velocity costates $\vec{p}_r(t_0)$ and $\vec{p}_v(t_0)$ are assessed by

$$\begin{cases} \vec{p}_v = p_v \vec{u} \\ \vec{p}_r = -\dot{\vec{p}}_v = -p_v \omega \vec{n} \end{cases} \tag{38}$$

The values needed for this assessment are $p_v$, $\vec{u}(t_0)$, $\vec{n}(t_0)$, $\omega(t_0)$.
The modulus of costate velocity is given by Eq. (15).

$$p_v = \frac{m_f^*}{v_e} \tag{39}$$

The projection of $\vec{u}$ and $\vec{n}$ in the Galilean frame (O,x,y) are (see Figure 1)

$$\vec{u} = \begin{pmatrix} \sin(\theta - \varphi) \\ \cos(\theta - \varphi) \end{pmatrix} \quad \text{and} \quad \vec{n} = \begin{pmatrix} -\cos(\theta - \varphi) \\ \sin(\theta - \varphi) \end{pmatrix} \tag{40}$$

The initial longitude $\varphi_0$ is given. The value of the initial pitch angle $\theta_0$ is obtained as the solution of Eq. (32) at the initial date. The angular rate $\omega_0$ in the Galilean frame is given by Eq. (25).

We obtain analytical formulae giving the initial costate components in the Galilean frame (O,x,y).

$$\vec{p}_v(t_0) = \frac{m_f^*}{v_e}\begin{pmatrix} \sin(\theta_0 - \varphi_0) \\ \cos(\theta_0 - \varphi_0) \end{pmatrix} \quad \text{and} \quad \vec{p}_r(t_0) = \omega_0 \frac{m_f^*}{v_e}\begin{pmatrix} -\cos(\theta_0 - \varphi_0) \\ \sin(\theta_0 - \varphi_0) \end{pmatrix} \tag{41}$$

We observe that the costates given by Eqs. (37,41) are proportional to the final mass $m_f^*$. Since the costate vectors introduced in the PMP equations are defined up to a multiplicative factor, knowing the optimal final mass is unnecessary for this assessment.

The initial costates $\vec{p}_r(t_0)$, $\vec{p}_v(t_0)$, $p_m(t_0)$ can therefore by assessed directly from the initial kinematic conditions $(r_0, v_0, \gamma_0, \varphi_0)$ <u>without solving any optimization problem</u>.

These initial costates can in turn be used as initial guess if the OCP Eq. (2) is to solve again with a fixed thrust level differing from the optimal one T*. In that case, the solution of Eq. (32) is no longer optimal. The costate equations Eq. (6) must be integrated to determine the optimal thrust direction by Eq. (5), and the initial costates become additional unknowns. This two boundary value problem (TBVP) is solved by a shooting method generally highly sensitive to the initial guess. The analytical guess derived from Eqs. (37,41) in the case of the optimal thrust level T* proves adequate to ensure the convergence of the shooting method in a large range of thrust levels around T*.

## 3. Application

The solution method is illustrated on a typical example of a launcher upper stage targeting a geostationary transfer orbit (GTO). The possible extensions to other similar problems are then discussed.

### 3.1 Illustrative Example

We consider a launcher upper stage with an initial gross mass of 40000 kg. The engine specific impulse is 350 s corresponding to an exhaust velocity of 3432.3 m/s. The initial conditions result from the flight of the launcher previous stage, which is injected on a fall-out trajectory. The apogee and perigee altitudes at the engine ignition are respectively 400 km and −5000 km and the initial anomaly is 169 deg. These orbital parameters correspond to an altitude of 164.5 km, a velocity of 4909.8 m/s and a flight path angle of 19.84 deg. The x axis of the reference Galilean frame is defined by the initial perigee. The initial longitude is thus equal to 169 deg.

The target is a geostationary transfer orbit (GTO) with an apogee at 36000 km and a perigee at 250 km. The Earth equatorial radius is $R_E$ = 6378137 m and the gravitational constant is $\mu = 3.986005.10^{14}$ m$^3$/s$^2$.

An initial guess for the thrust level T and the final time $t_f$ is derived using the Eqs. (34−36). The trajectory is simulated for this initial guess and using the closed-loop command law resulting from Eq. (32). The third column of Table 1 presents these starting values, and the apogee and perigee reached at the end of the trajectory. The nonlinear system with 2 unknowns (T and $t_f$) and 2 equations (apogee and perigee) is easily solved by a Newton–like method. The solution obtained is given in the last column of Table 1. This solution is quite close to the initial guess favoring the fast convergence of the Newton-like method.

The initial costates of the OCP are retrieved analytically from Eqs. (37,41). The initial pitch angle $\theta_0$ is found by solving Eq. (32) at the initial date. It depends only on the kinematic conditions $r_0$, $v_0$, $\gamma_0$. The initial angular rate $\omega_0$ is assessed from Eq. (25).

$$\begin{cases} r_0 = 6542.632\,\text{km} \\ v_0 = 4909.82\,\text{m/s} \\ \gamma_0 = 19.8414\,\text{deg} \end{cases} \Rightarrow \theta_0 = 7.5153\,\text{deg} \;,\; \omega_0 = -0.06658\,\text{deg/s} \tag{42}$$

Table 1 gives the costate components in the Galilean frame. The difference between the third and the last column comes only from the multiplicative factor $m_f$.

| Variable | Unit | Initial guess | Optimal solution |
|---|---|---|---|
| Thrust level | kN | 126.038 | 131.296 |
| Final time | s | 856.0 | 839.19 |
| Final mass | kg | 8568.2 | 7898.4 |
| Apogee altitude | km | 24026.2 | 36000.0 |
| Perigee altitude | km | 188.4 | 250.0 |
| Position costate $p_x, p_y$ | kg/km | (−2.7502, 0.9210) | (−2.5354, 0.8491) |
| Velocity costate $p_{vx}, p_{vy}$ | kg/(m/s) | (−0.7927, −2.3670) | (−0.7308, −2.1821) |
| Mass costate $p_m$ | kg/kg | 0.2142 | 0.1975 |

Table 1 : Initial guess and optimal solution

The Figure 3 depicts the evolution of the altitude, the velocity, the flight path angle, the local pitch angle, the osculating apogee and perigee along the optimal trajectory.

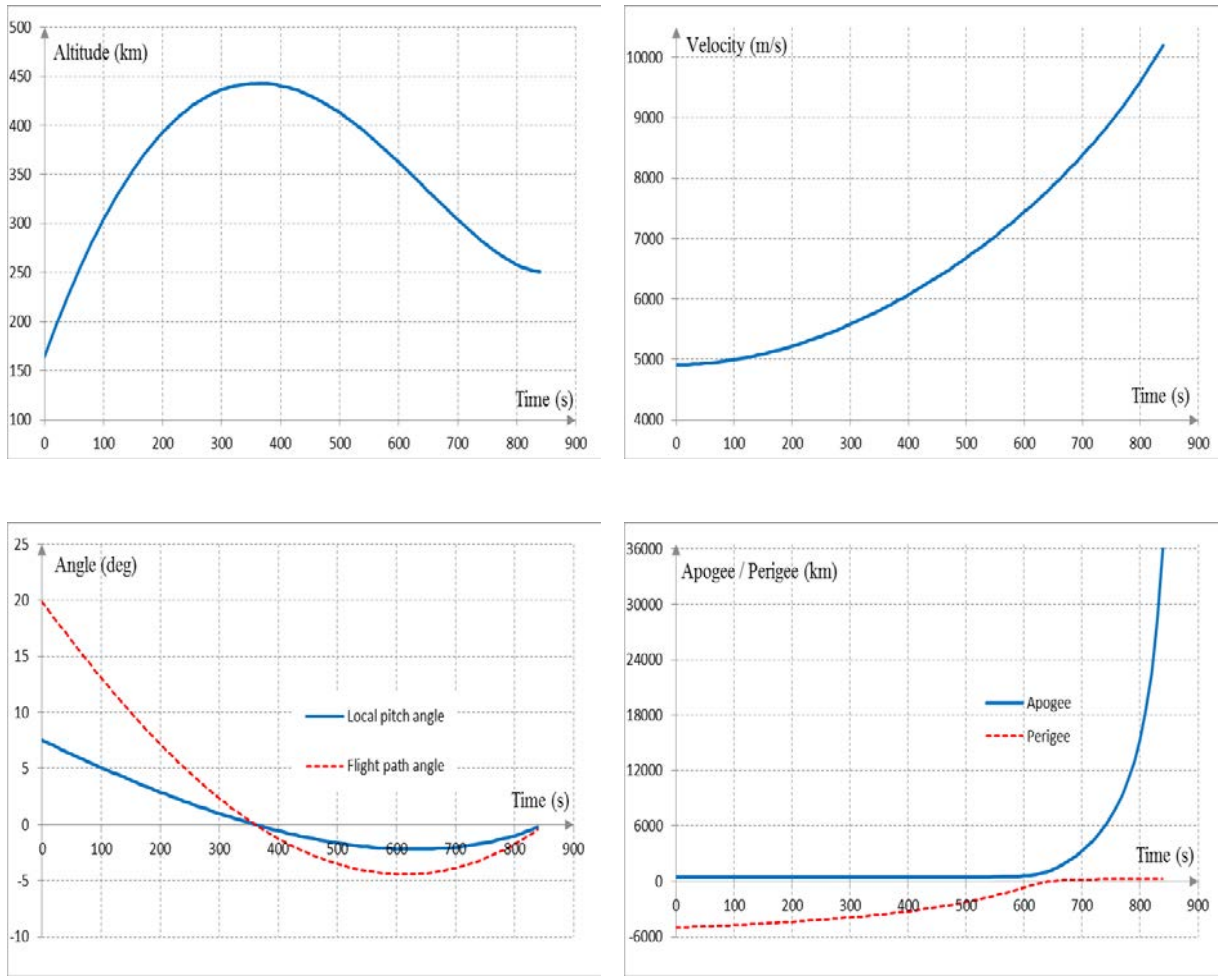

Figure 3 : Optimal trajectory variables

## 3.2 Sensitivity to the Thrust Level

The above example is solved again for thrust levels varying in the range 100 kN to 230 kN. For a fixed thrust level differing from the optimal one, the closed-loop control given by Eq. (32) is no longer optimal. It is necessary to integrate the costate equations Eq. (6) to retrieve the optimal thrust direction from Eq. (5) and therefore to solve the shooting problem with unknown initial costate values and with endpoints conditions given by Eq. (7). Finding a good initial guess for this problem is usually very challenging and it is the major drawback of the indirect approach. This drawback can be overcome owing to the analytical costates given by Eqs. (37,41) corresponding to the optimal thrust level. These costates are assessed up to a multiplicative factor corresponding to the final mass, and they depend only on the initial kinematic conditions. They are thus available without solving any optimization problem. This starting solution allows solving easily a series of shooting problems by a homotopic process on the successive thrust levels in the desired range.

The Figure 4 depicts the final mass and the velocity losses depending on the thrust level. The best performance is effectively recovered for the thrust level of 131 kN previously found in §3.1.

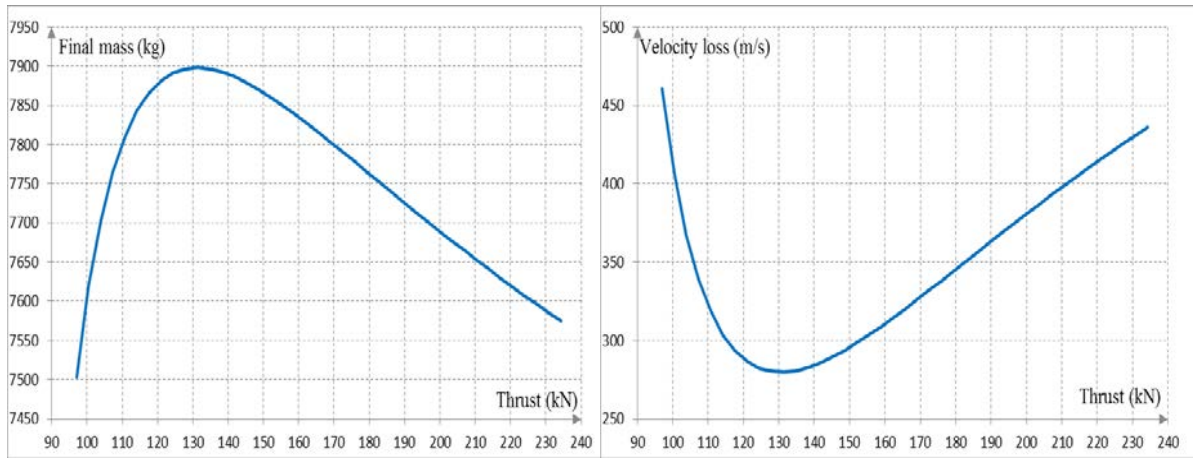

Figure 4 : Performance vs thrust level

From Eq. (9) the Hamiltonian of the OCP Eq. (2) is constantly null whatever the thrust level. The Figure 5 depicts the components $H_0$ and $T.H_T$ of the Hamiltonian written as $H = H_0 + T.H_T$ in Eq. (23). These components are plotted at the initial date for thrust levels varying from 100 kN to 230 kN. We observe as expected that their sum is null for every thrust level, and that both components vanish for the optimal thrust level.

The initial costate components $p_m$ (black bold line), $p_{rx}$ and $p_{ry}$ (blue thin solid and dotted lines), $p_{vx}$ and $p_{vy}$ (red medium solid and dotted lines) are also plotted on the Figure 5. The respective units are kg/kg for $p_m$, kg/km for $p_r$, and kg/(m/s) for $p_v$. The smooth evolution of the initial costates explains the good convergence of the shooting method when applying a homotopic process on the thrust level.

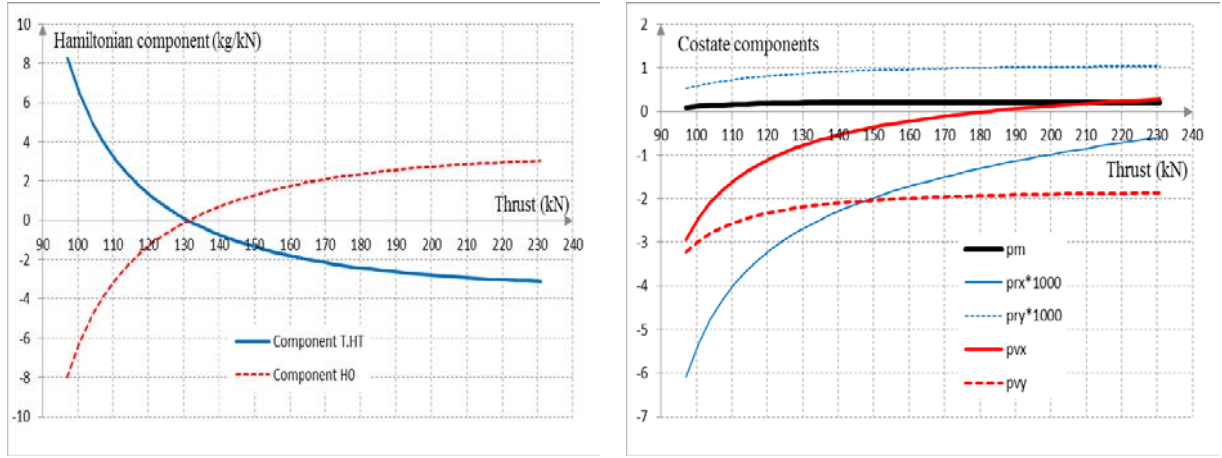

Figure 5 : Initial Hamiltonian and costate components vs thrust level

### 3.3 Extensions

The solution method presented in section §3.1 applies to a launcher upper stage planar trajectory with a constant optimized thrust level. The analysis can be extended to similar minimum-fuel problems by relaxing some of these assumptions.

A first extension is to consider that the thrust law is not necessarily constant, but that it depends both on the time since the ignition and on a scalar parameter $\tau$ that can be optimized. With a thrust law of the form $T(t-t_0,\tau)$, the condition Eq. (10) becomes

$$\int_{t_0}^{t_f} \frac{\partial H}{\partial \tau} dt = 0 \quad \Rightarrow \quad \int_{t_0}^{t_f} \frac{\partial H}{\partial T} \frac{\partial T}{\partial \tau} dt = 0 \quad \Rightarrow \quad \int_{t_0}^{t_f} \left( \frac{p_v}{m} - \frac{p_m}{v_e} \right) \frac{\partial T}{\partial \tau} dt = 0 \tag{43}$$

If the derivative of T with respect to $\tau$ does not vanish identically on a time interval, we recover the condition Eq. (12) and the subsequent analysis holds unchanged. The solution method can thus be applied for example when the thrust law follows a predefined shape $T_0(t-t_0)$ with an additive or a multiplicative variation parametrized by $\tau$.

A second extension is to consider non planar trajectories. The upper stage of a launcher is generally injected directly in the targeted orbital plane, either by a proper choice of the launch azimuth, or by out-of-plane thrusting of the previous stages. Nevertheless some out-of-plane thrusting may be still necessary during the upper stage flight, for example when there is a constraint on the perigee argument. The analysis made in the planar case can be followed until Eq. (19) and the constant modulus property of the costate velocity still holds. At this step it becomes necessary to account for an out-of-plane angle $\psi$ between the thrust direction $\vec{u}$ and the trajectory plane $(\vec{r},\vec{v})$. This unknown angle will be present in the Eqs. (21,25) and the closed-loop in-plane command will depend on $\psi$. In practice, the out-of-plane thrust component (if any) is quite small during an upper stage flight. Considering a constant angle $\psi$ as additional optimization parameter allows still using the closed-loop solution for the in-plane command without recourse to the costate equations. The solution method can thus be applied to non-planar transfers with an inclination constraint.

The third extension concerns low-thrust transfers. The analysis applies without change and the condition Eq. (30) satisfied by the optimal thrust direction is still valid. This equation has at least 2 solutions yielding opposite thrust directions, and it may even have 4 solutions. For a launcher upper stage with a high thrust engine, the first solution corresponding to Eq. (32) is obviously the optimal one. This is no longer true for low-thrust transfers which can present thrusting arcs opposite to the velocity. Unfortunately the PMP only provides local necessary conditions and there is no rigourous mean to determine locally which solution is the optimal one among the 2 or 4 candidates. For such transfers the correct thrust direction can be recovered from Eq. (5) by integrating the $\vec{p}_r$ and $\vec{p}_v$ costate equations Eq. (16). Their initial values are given analytically by Eq. (41) up to the multiplicative factor $m_f$, but this factor does not matter since the costate equations are homogeneous and since the thrust direction $\vec{u}$ does not depend on the costate norm. The initial guess issue is thus avoided and the optimal low-thrust trajectory can be obtained as previously by solving a nonlinear system with only the thrust level and the final time as unknowns.

## 4. Conclusion

The minimum-fuel trajectory of a launcher upper stage has been investigated in the case of a constantly thrusting engine. The thrust level is constant and its value has to be optimized together with the thrust direction.

This optimization parameter induces an additional optimality condition to the PMP equations. Exploiting this condition, a closed-loop control can be derived in the planar case for the thrust direction. Guessing the initial costate and propagating the costate equations is no longer necessary. The problem is reduced to a nonlinear system with the thrust level and the final time as unknowns. A rough initial estimate of these two parameters is built from physical considerations on the required velocity impulse and on the vehicle acceleration level. The optimal trajectory and the associated thrust level are thus found very easily.

On the other hand the initial costate components are assessed analytically from the initial conditions. These costates can in turn be used as initial guess when the minimum-fuel problem must be solved again for different thrust levels, differing from the optimal one. For these problems, the closed-loop solution is no longer valid and a shooting method must be applied. The analytical costates obtained for the optimal thrust level prove a quite adequate guess to ensure the fast convergence of the shooting method.

The closed-loop solution proposed in this paper allows solving easily the minimum-fuel launcher trajectory in the planar case, when the engine thrust level is constant. The method can similarly be applied in the case of any thrust law defined by one parameter. The extension to non planar trajectories can also be envisioned by adding an out-of-plane parametric command to the optimization parameters, while keeping the in-plane closed-loop solution. Minimum-fuel low-thrust transfers can also be adressed by the same solution method, with the difference that the closed-loop solution is no longer fully determined from the PMP. Nevertheless the optimal thrust direction can still be recovered by the costate propagation starting from the analytical initial values. The initial guess issue is thus avoided and the optimal control problem is again reduced to a 2-unknowns nonlinear system.

**References**

[1] Goddard, R.H., "A method of reaching extreme altitudes", Smithsonian Misc. Collect. **71**(2), Smithsonian institution, Washington (1919).
[2] Bryson, A.E., Ho, Y.C. *Applied optimal control*, Hemisphere Publishing Corporation, 1975.
[3] Lawden, D.F. *Optimal trajectories for space navigation*, Butterworths Publishing Corporation, 1963.
[4] Leitmann, G., *Optimization technique with applications to Aerospace System*, Academic Press New-York, 1962.
[5] Cerf, M., Haberkorn, T., Trélat, E., "Continuation from a flat to a round Earth model in the coplanar orbit transfer problem", *Optimal Control and Applied Mathematics*, DOI: 10.1002/oca.1016, 2011.
[6] Betts, J.T., *Practical methods for optimal control and estimation using nonlinear programming*, Siam, 2010.
[7] Conway, B.A., *Spacecraft Trajectory Optimization*, Cambridge University Press, 2010.
[8] P. Martinon, J. Gergaud, "Using switching detection and variational equations for the shooting method", *Optimal Cont. Appl. Methods* 28, no. 2 (2007), 95--116.
[9] Bonnans, F., Martinon, P., Trélat, E., "Singular arcs in the generalized Goddard's problem", *Journal of Optimization Theory and Applications*, Vol 139 (2), pp 439-461, 2008
[10] Pontryagin, L., Boltyanskii, V., Gramkrelidze, R., Mischenko, E., *The mathematical theory of optimal processes*, Wiley Interscience, 1962.
[11] Trélat, E., *Contrôle optimal – Théorie et Applications*, Vuibert, 2005.
[12] Hull, D.G., *Optimal control theory for applications*, Springer, 2003.
[13] Boltyanski, V.G., Poznyak, A.S., *The robust maximum principle*, Birkhäuser, 2012.
[14] Chobotov, V., *Orbital mechanics*, 3rd Edition, AIAA Education Series, 2002
[15] Tewari, A., *Advanced control of spacecraft and aircrafts*, Siam, 2010.

**Acronyms**

| | |
|---|---|
| OCP | Optimal Control Problem |
| PMP | Pontryaguin Maximum Principle |
| TPBVP | Two Point Boundary Value Problem |
| NLP | Non Linear Programming |
| GTO | Geostationary Transfer Orbit |